\documentclass{amsart}
\usepackage{amssymb}
\usepackage{verbatim}
\usepackage{array}
\usepackage{latexsym}
\usepackage{enumerate}
\usepackage{amsmath}
\usepackage{amsfonts}
\usepackage{graphicx,epsfig}
\usepackage[english]{babel}
\newcommand{\PG}{\rm {PG}}
\newcommand{\K}{\mathbb{K}}

\begin{document}

\title{3-nets realizing a diassociative loop in a projective plane}

\author{G\'abor Korchm\'aros}
\address{Dipartimento di Matematica, Informatica ed Economia\\
Universit\`a della Basilicata\\
Contrada Macchia Romana\\
85100 Potenza, Italy}
\email{gabor.korchmaros@unibas.it}

\author{G\'abor P. Nagy}
\address{Bolyai Institute \\
University of Szeged \\
Aradi v\'ertan\'uk tere 1\\
H-6720 Szeged, Hungary}
\address{MTA-ELTE Geometric and Algebraic Combinatorics Research Group \\
P\'azm\'any P. s\'et\'any 1/c \\
H-1117 Budapest, Hungary}
\email{nagyg@math.u-szeged.hu}

\newtheorem{theorem}{Theorem}[section]
\newtheorem{proposition}[theorem]{Proposition}
\newtheorem{lemma}[theorem]{Lemma}
\newtheorem{corollary}[theorem]{Corollary}
\newtheorem{scholium}[theorem]{Scholium}
\newtheorem{observation}[theorem]{Observation}
\newtheorem{assertion}[theorem]{Assertion}
\newtheorem{result}[theorem]{Result}
{\theoremstyle{definition}
\newtheorem*{definition*}{Definition}
\newtheorem{example}[theorem]{Example}
\newtheorem{remark}[theorem]{Remark}
\newtheorem*{proposition*}{Proposition}
\newtheorem*{corollary*}{Corollary}
\newtheorem*{lemma*}{Lemma}

\begin{abstract}
A \emph{$3$-net} of order $n$ is a finite incidence structure consisting of points and three pairwise disjoint classes of lines, each of size $n$, such that every point incident with two lines from distinct classes is incident with exactly one line from each of the three classes. The current interest around  $3$-nets (embedded) in a projective plane $\PG(2,\K)$, defined over a field $\K$ of characteristic $p$, arose from algebraic geometry; see  \cite{fy2007,miq,per,ys2004,ys2007}. It is not difficult to find $3$-nets in $\PG(2,\K)$ as far as $0<p\le n$. However, only a few infinite families of $3$-nets in $PG(2,\K)$ are known to exist whenever  $p=0$, or $p>n$. Under this condition, the known families are characterized as the only $3$-nets in $\PG(2,\K)$ which can be coordinatized by a group; see \cite{knp_3}. In this paper we deal with $3$-nets in $PG(2,\K)$ which can be coordinatized by a diassociative loop $G$ but not by a group. We prove two structural theorems on $G$. As a corollary, if $G$ is commutative then every non-trivial element of $G$ has the same order, and $G$ has exponent $2$ or $3$. We also discuss the existence problem for such $3$-nets.
\end{abstract}

\maketitle


\section{Introduction}
The concept of a $3$-net comes from classical differential geometry via the combinatorial abstraction of the concept of a $3$-web. Formally, a \emph{$3$-net} of order $n$ is a finite incidence structure consisting of points and three pairwise disjoint classes of lines, each of size $n$, such that every point incident with two lines from distinct classes is incident with exactly one line from each of the three classes. It is well known that every $3$-net can be coordinatized by a loop. The set $Q$ endowed with a binary operation ``$\cdot$'' is a \textit{quasigroup,} if for any $a,b \in Q$, the equations $a\cdot x=b$ and $y\cdot a=b$ have unique solutions in $Q$. A quasigroup with a multiplicative unit element is called a \textit{loop.} For a general reference on nets, loops and quasigroups see for instance \cite{Barl,bruck}.

In this paper we deal with $3$-nets (embedded) in $\PG(2,\K)$, the projective plane over a field $\K$ of characteristic $p\geq 0$. Such a $3$-net, with line classes ${\mathcal{A}},{\mathcal{B}},{\mathcal{C}}$ and coordinatizing loop $G=(G,\cdot)$, is equivalently defined by a triple of bijective maps from $G$ to $({\mathcal{A}},{\mathcal{B}},{\mathcal{C}})$, say
$$\alpha:\,G\to {\mathcal{A}},\,\beta:\,G\to {\mathcal{B}},\,\gamma:\,G\to {\mathcal{C}}$$
such that $a\cdot b=c$ if and only if $\alpha(a),\beta(b),\gamma(c)$ are three concurrent lines in $PG(2,\mathbb{K})$, for any $a,b,c \in G$. If this is the case, the $3$-net in $PG(2,\mathbb{K})$ is said to {\em{realize}} the loop $G$.

For the purpose of investigating $3$-nets in $\PG(2,\K)$, the groundfield $\K$ may be assumed to be algebraically closed. In order to present the key examples of embedded $3$-nets, it is convenient to work with the dual concept. Formally, a \textit{dual $3$-net} of order $n$ in $PG(2,\mathbb{K})$ consists of a triple $(\Lambda_1,\Lambda_2,\Lambda_3)$ with $\Lambda_1,\Lambda_2,\Lambda_3$ pairwise disjoint point-sets of size $n$, called {\em{components}}, such that every line meeting two distinct components meets each component in precisely one point. We notice that finite dual $3$-nets are also called \textit{transversal designs.}

The following concepts and results have a detailed exposition in \cite{knp_3}. We say that an embedded dual $3$-net is \textit{algebraic,} if its point set $\Lambda_1 \cup \Lambda_2 \cup \Lambda_3$ is contained in a cubic curve $\mathcal{F}$. If $\mathcal{F}$ is reducible then we speak of \textit{pencil type,} \textit{triangular type} or \textit{conic-line type} dual $3$-net. Except for the pencil type, all algebraic (dual) $3$-nets are coordinatized by either a cyclic group or by a direct product of two cyclic groups. Finite dihedral groups can be realized by dual $3$-nets of \textit{tetrahedron type;} in this case the point set is contained in six lines joining four independent points. Finally, we mention that the quaternion group $\mathbf{Q}_8$ has an exceptional realization, cf. \cite{urzua2009}.

In recent years, finite $3$-nets realizing a group have been investigated also in connection with complex line arrangements and resonance theory; see \cite{bkm,bkn,fy2007,knp_3,knp_k,miq,per,ys2004,ys2007}. The following almost complete classification of such $3$-nets is proven in \cite{knp_3}.
\begin{theorem}
\label{mainteo} In the projective plane $PG(2,\mathbb{K})$ defined over an algebraically closed field $\mathbb{K}$ of characteristic $p\geq 0$, let $(\Lambda_1,\Lambda_2,\Lambda_3)$ be a dual $3$-net of order $n\geq 4$ which realizes a group $G$. If either $p=0$ or $p>n$ then one of the following holds.
\begin{itemize}
\item[\rm{(I)}] $G$ is either cyclic or the direct product of two cyclic groups, and $(\Lambda_1, \Lambda_2, \Lambda_3)$ is algebraic.
\item[\rm{(II)}] $G$ is dihedral and $(\Lambda_1,\Lambda_2,\Lambda_3)$ is of tetrahedron type.
\item[\rm{(III)}] $G$ is the quaternion group of order $8$.
\item[\rm{(IV)}] $G$ has order $12$ and is isomorphic to $\rm{Alt}_4$.
\item[\rm{(V)}] $G$ has order $24$ and is isomorphic to $\rm{Sym}_4$.
\item[\rm{(VI)}] $G$ has order $60$ and is isomorphic to $\rm{Alt}_5$.
\end{itemize}
\end{theorem}
A computer aided exhaustive search shows that if $p=0$ then (IV) (and hence (V), (VI)) does not occur; see \cite{np}. It has been conjectured that this holds true in any characteristic.

In this paper we focus on  $3$-nets in $PG(2,\K)$ which can be coordinatized by a diassociative loop $G$ different from a group. Recall that a loop $G$ is \textit{diassociative} if any subloop generated by two elements is a group. There are two important classes of diassociative loops: Moufang loops and Steiner loops. Moufang loops are loops satisfying one (hence all) of the following identities.
\[ z(x(zy)) = ((zx)z)y, \hspace{1cm} x(z(yz)) = ((xz)y)z, \hspace{1cm} (zx)(yz) = (z(xy))z.\]
In general, Moufang loops have a rich algebraic structure. This is not the case for \textit{Steiner loops.} Steiner loops are diassociative loops of exponent two. Finite Steiner loops are in one-to-one connection with Steiner triple systems. For other classes of diassociative loops we refer to \cite{KinyonKunnenPhillips}.

Our results consist of three  structural theorems on $G$.
\begin{theorem}
\label{prop1} In the projective plane $PG(2,\mathbb{K})$ defined over an algebraically closed field $\mathbb{K}$ of characteristic $p\geq 0$, let $(\Lambda_1,\Lambda_2,\Lambda_3)$ be a dual $3$-net of order $n\geq 4$ which realizes a diassociative loop $G$ different from a group. Let $d$ be the
maximum of the orders of the elements in $G$, and suppose that $d\geq 4$. If either $p=0$ or $p>n$ then one of the following holds.
\begin{itemize}
\item[\rm(a)] $G$ has a unique subgroup $H$ of order $d$. Moreover, each element not in $H$ is an involution, and two such involutions either commute or their product is in $H$.
\item[\rm(b)] $d=4$, and $G$ has a subgroup isomorphic to one of the groups ${\bf{Q}}_8,\,\rm{Alt}_4$.
\end{itemize}
\end{theorem}
\begin{theorem}
\label{prop2} With the same hypotheses as in Theorem {\emph{\ref{prop1}}}, assume further that $G$ contains a subgroup isomorphic to  ${\bf{Q}}_8$ but no subgroup isomorphic to $\rm{Alt}_4$.
Then $G$ has a unique involution and the subgroup generated by any two non-commuting elements is isomorphic to ${\bf{Q}}_8$.
\end{theorem}
It may be observed that a loop $G$ as in Theorem \ref{prop2} defines a Steiner triple system in a natural way, namely the points are subgroups of $G$ of order $4$ and the blocks are the subgroups isomorphic to ${\bf{Q}}_8$, the point-block incidence being the set theoretic inclusion.

For a commutative loop $G$, neither (a) nor (b) of Theorem \ref{prop1} can occur, and hence $d\leq 3$. More precisely, the following result holds.
\begin{corollary}
\label{cor1}  With the same hypotheses as in Theorem {\emph{\ref{prop1}}}, assume further that $G$ is commutative. Then every non-trivial element in $G$ has the same order, and $G$ has exponent $2$ or $3$.
\end{corollary}
The quaternion group $\mathbf{Q}_8$ has a counterpart in the class of Moufang loops. Let $\mathbb{O}$ be the division ring of real octonions and let $1,e_1,\ldots,e_7$ be an orthonormal basis. The set
\[\mathbf{O}_{16}=\{\pm 1,\pm e_1,\ldots,\pm e_7 \}\]
forms a Moufang loop with a unique involution $-1$ and $14$ elements of order $4$. ($\mathbf{O}_{16}$ is also called the \textit{Cayley loop of order 16}.)
\begin{theorem}
\label{mo}With the same hypotheses as in Theorem {\emph{\ref{prop1}}}, assume further that $G$ is a Moufang loop. Then $G$ contains either the octonion loop $\mathbf{O}_{16}$, or it has a subgroup isomorphic to $\rm{Alt}_4$.
\end{theorem}
An interesting issue which appears to be rather difficult is the existence and construction of $3$-nets in the classical projective plane $\PG(2,\K)$ realizing a loop different from a group. All such examples available in the literature are $3$-nets of order $n=5,6$, obtained by computer aided searches; see \cite{urzua2009}.

\section{Proof of Theorem \ref{prop1}}
Let $\Lambda=(\Lambda_1,\Lambda_2,\Lambda_3)$ be a $3$-net of order $n$ coordinatized by a diassociative loop $G$ but not by a group.  Let $g\in G$ be an element whose order is $d$, and let $\Psi=(\Psi_1,\Psi_2,\Psi_3)$ be the $3$-net (subnet of $\Lambda$) coordinatized by $\langle g \rangle$. Take any element $h\in G$ not lying in the cyclic group generated by $g$, and consider the subgroup $H$ generated by $g$ and $h$. Obviously, $H$ is not a cyclic group. Since $H$ is a subloop of $G$, $H$ also realizes a $3$-net $\Delta=(\Delta_1,\Delta_2,\Delta_3)$ in $PG(2,\K)$. The classification \cite[Theorem 1.1]{knp_3} applies to $H$ and yields one of the cases below, apart from the sporadic cases. Therefore, dismissing (b), we have that either
\begin{itemize}
\item[(i)] $H$ is the direct product of two cyclic subgroups, and $\Delta_1\cup\Delta_2\cup \Delta_3$ lies on a plane non-singular cubic curve $\mathcal{F}_3$, or
\item[(ii)] $H$ is a dihedral group, and $\Delta$ is of tetrahedron type.
\end{itemize}
We first investigate case (i). As $G$ is not a group, it must contain an element $u\not\in H$. Replacing $h,H$ by $u,U=\langle g,u \rangle$ in the above argument shows that $U$ realizes a $3$-net $\Phi=(\Phi_1\cup \Phi_2\cup \Phi_3)$ in $PG(2,\K)$, and that $U$ is either the direct product of two cyclic groups, or it is a dihedral group. The latter case here cannot actually occur. To show this, observe that if $U$ is dihedral then the maximality of $d$ implies that $u$ is an (involutory) element lying in some coset of $\langle g \rangle$. Since $\Phi$ is of tetrahedron type in this case,
we have that $\Psi$ is a triangular $3$-net in $PG(2,\K)$ of order $d$. But this is impossible in our case, since the points of $\Psi$ lie on $\mathcal{F}_3$. In fact, $\mathcal{F}_3$ is non-singular while $\Psi_1$ consists of $d>3$ collinear points.
Therefore, $U$ is a direct product of two cyclic groups and $\Phi_1\cup\Phi_2\cup \Phi_3$ lies on a non-singular plane cubic curve $\mathcal{F}_1$. The intersection $\mathcal{F}_3\cap \mathcal{F}_1$ contains all points of $\Psi$. Since $d>3$, this yields $\mathcal{F}_3=\mathcal{F}_1$. Since $u$ denotes any element of $G$ outside $H$, it turns out that $\Lambda_1\cup\Lambda_2\cup\Lambda_3$ lies on $\mathcal{F}_3$. But then $G$ itself is a group, the direct product of two cyclic groups. Therefore, (i) cannot actually occur.

In case (ii), we may assume that $H_1=\langle g,h_1 \rangle$ is dihedral for any $h_1\in G\setminus\langle g \rangle$. Therefore, $h_1$ is an involution. Moreover, if $h_2$ is another involution in $G$, then either $h_2$ commutes with $h_1$, or $h_2$ lies $H_1$.

\section{Proof of Theorem \ref{prop2}}
From the definition of a dual $3$-net, there is a triple of bijective maps from $G$ to $(\Lambda_1,\Lambda_2,\Lambda_3)$, say $\alpha,\beta,\gamma$ respectively such that $a\cdot b=c$ in $G$ if and only if $\alpha(a),\beta(b),\gamma(c)$ are three collinear points in $PG(2,\mathbb{K})$, for any $a,b,c \in G$.

Choose two elements $g_1,g_2\in G$ which generate a subgroup $H$ isomorphic to ${\bf{Q}}_8$. We remark that case (a) in Theorem \ref{prop1} cannot occur since both $g$ and $h$ have order $4$. Therefore, $d=4$. Set $g_3=g_1g_2$; then $\langle g_1 \rangle$, $\langle g_2 \rangle$, $\langle g_3 \rangle$ are the three cyclic subgroups of order $4$ in $H$. Take an element $u\in G$ not lying in $H$.


Assume that $u$ is an involution. For $i=1,2,3$, the group $U_i=\langle u,g_i \rangle$ contains at least two distinct involutions, and hence it is not isomorphic to ${\bf{Q}}_8$. From Theorem \ref{mainteo} applied to $U_i$, we deduce that either $U_i$ is dihedral, or abelian.

We first investigate the case when all $U_i$'s are abelian. Clearly, $U_i=\langle g_i \rangle \times \langle u \rangle \cong C_4\times C_2$, and case (a) of Theorem \ref{prop1} holds for the sub $3$-net $(\Delta_1^i,\Delta_2^i,\Delta_3^i)$ realizing $U_i$. Let $\mathcal{F}_i$ be the cubic curve containing $\Delta_1^i\cup\Delta_2^i\cup\Delta_3^i$. Since $U_i$ is not cyclic, $\mathcal{F}_i$ is nonsingular. These three sub $3$-nets of order $8$ share a sub $3$-net of order $4$, say $(\Omega_1,\Omega_2,\Omega_3)$, realizing the group $T=\langle u,g_1^2=g_2^2\rangle$. Since $|\mathcal{F}_1\cap \mathcal{F}_2\cap \mathcal{F}_3|\geq 12>9$, this yields
that $\mathcal{F}_1=\mathcal{F}_2=\mathcal{F}_3$. But then the sub $3$-net realizing $H$ lies on $\mathcal{F}_1$ a contradiction since ${\bf{Q}}_8$ is not abelian.

Assume now that $U_1,U_2$ are abelian and $U_3$ dihedral. With the same argument, the dual sub $3$-nets realizing $U_1,U_2$ are contained in the nonsingular cubic curve $\mathcal{F}$. The dual sub $3$-net realizing $U_3$ is of tetrahedron type, which means that the four points of $\alpha(\langle g_3 \rangle)$ are collinear. The triple $(\alpha(\langle g_3 \rangle), \beta(H\setminus \langle g_3 \rangle), \gamma(H\setminus \langle g_3 \rangle))$ is a dual $3$-net realizing $\langle g_3 \rangle$. On the one hand, the 8 points of $\beta(H\setminus \langle g_3 \rangle) \cup \gamma(H\setminus \langle g_3 \rangle)$ are contained in a (possibly degenerate) conic $\mathcal{C}$, see \cite[Theorem 5.1]{bkm}. On the other hand, $H\setminus \langle g_3 \rangle \subset \langle g_1 \rangle \cup \langle g_2 \rangle \subset U_1\cup U_2$. This implies $\beta(H\setminus \langle g_3 \rangle) \cup \gamma(H\setminus \langle g_3 \rangle) \subset \mathcal{F}$ and $|\mathcal{F} \cap \mathcal{C}|\geq 8$, a contradiction.

Assume that $U_1,U_2$ are dihedral. Hence the dual sub $3$-nets realizing $U_1,U_2$ are of tetrahedron type yielding that the four points of $\alpha(\langle g_1 \rangle)$ and the four points of $\alpha(\langle g_2 \rangle)$ are contained in the lines $\ell_1,\ell_2$, respectively. However, $\alpha(1), \alpha({g_1^2=g_2^2}) \in \ell_1 \cap \ell_2$, thus, $\ell_1=\ell_2$. Similarly, the six points of $\beta(\langle g_1 \rangle \cup \langle g_2 \rangle)$ and the six points of $\gamma(\langle g_1 \rangle \cup \langle g_2 \rangle)$ are contained in the lines $m,m'$, respectively. If $U_3$ is dihedral, then the dual sub $3$-net realizing $H$ is contained in $\ell_1 \cup m \cup m'$, which is impossible since $H$ is not cyclic. If $U_3$ is abelian, then the sub $3$-net realizing it is contained in the nonsingular cubic curve $\mathcal{F}$. The second component of its sub $3$-net $(\alpha(\langle g_3 \rangle), \beta(H\setminus \langle g_3 \rangle), \gamma(H\setminus \langle g_3 \rangle))$ is contained in $m$, hence $|m\cap \mathcal{F}|\geq 4$, a contradiction.


Assume that $u$ has order $3$. Then $U$ is neither a dihedral group nor isomorphic to ${\bf{Q}}_8$. From Theorem \ref{mainteo},  $U=\langle g\rangle \times \langle u \rangle$ and hence $U$ is a cyclic group of order $12$ contradicting the remark at the beginning about case (a) in Theorem \ref{prop1}.

Therefore, $G$ contains just one involution $v$, and if $u\neq v$ then $u$ has order $4$. Let $u_1,u_2\in G$ be any two distinct elements other than $v$. Since $U$ contains no element of order $3$, $U=\langle u_1,u_2 \rangle$ is a $2$-group of exponent $4$ containing a unique involution. Since $U$ has order bigger than $4$, the only possibility is $U\cong  {\bf{Q}}_8$.

\begin{remark}
Let $S$ be the Steiner loop of order $10$ corresponding to the Steiner triple system $AG(2,3)$. $S$ has a central extension $Q$ of order $20$ all proper subloops are isomorphic to $C_2,C_4,C_2\times C_4$, or ${\bf{Q}}_8$. In particular, $Q$ is diassociative. By Theorem \ref{prop2}, $Q$ has no projective realization despite all its subloops have.
\end{remark}


\section{Proof of Theorem \ref{mo}}

We start with three important facts on Moufang loops of small exponent. First, as diassociative loops of exponent $2$ are commutative, Moufang loops of exponent $2$ are elementary abelian groups. Second, by \cite[Corollary 1]{Glauberman} finite Moufang loops of exponent $3$ are nilpotent. This implies that any proper finite Moufang loop of exponent $3$ contains a subloop of order $27$. The classification of small Moufang loops \cite{Go} shows that Moufang loops of order $27$ are groups. Thus, if $G$ is a Moufang loop of exponent $3$ then it contains a subgroup $H$ of order $27$. Since no such group $H$ has a realization in $\PG(2,\K)$ by Theorem \ref{mainteo}, we have a contradiction.

Let us assume that $G$ has an element $g$ of order $d>4$. Put $U=\langle g \rangle$. By Theorem \ref{prop1}, any $h \in G\setminus U$ has order $2$ and $\langle U,h \rangle$ is a dihedral group of order $2d$. In particular, on the one hand, $hU = Uh$, and on the other hand, the involutions generate $G$. \cite[Theorem 1]{Gagola} implies that $U$ is a normal subloop of $G$. For any subset $X$ of $G$, let $\Lambda_i(X)$ denote the points of $\Lambda_i$, indexed by the elements of $X$. \cite[Proposition 22]{knp_3} implies that any of the sets $\Lambda_i(U)$, $\Lambda_i(Uh)$ is contained in a line.

Choose an element $h \in G\setminus U$. Then we have four points $P,Q,R,S$ such that the point sets $\Lambda_1(U), \Lambda_2(U), \Lambda_3(U), \Lambda_1(Uh), \Lambda_2(Uh),\Lambda_3(Uh)$ are contained in the lines $QR, RS, PR, SP,SQ,PQ$. In fact, the points $P,Q,R,S$ are the vertices of the tetrahedron type dual $3$-net, corresponding to the dihedral group $\langle U,h \rangle$. Only the vertex $S$ depends on the choice of $h$; $S=S_h$.

Choose elements $h_1,h_2 \in G\setminus U$ such that $\langle U,h_1,h_2\rangle$ is a non-associative subloop of $G$. Let $P,Q,R,S_{h_1},S_{h_2},S_{h_1h_2}$ be the vertices of the tetrahedron type dual nets of $\langle U,h_1\rangle$, $\langle U,h_2\rangle$ and $\langle U,h_1h_2\rangle$. The sets
\[\Lambda_1(Uh_1), \Lambda_2(Uh_2), \Lambda_3(Uh_1h_2)\]
of points form a triangular dual 3-net. \cite[Proposition 10]{knp_3} implies $S_{h_1}=S_{h_2}=S_{h_1h_2}$, a contradiction.

Finally, assume that $G$ has no subgroup isomorphic to $\rm{Alt}_4$. By Theorem \ref{prop2}, $G$ has a unique (central) involution $u$. As two non-commuting elements generate a subloop isomorphic to $\mathbf{Q}_8$, the factor $G/\langle u\rangle$ is an elementary abelian $2$-group. Thus, $G$ contains a non-associative subloop $S$ of order $16$ with a unique involution. Using the classification of small Moufang loops in \cite{Go}, $S\cong \mathbf{O}_{16}$ follows.

\section{acknowledgement}
The work has been carried out within the Project PRIN (MIUR, Italy) and GNSAGA. The publication is supported by the European Union and co-funded by the European Social Fund. Project title: \textit{Telemedicine-focused research activities on the field of Mathematics, Informatics and Medical sciences.} Project number: TAMOP-4.2.2.A-11/1/KONV-2012-0073.

\end{document}